\documentclass[a4paper,11pt]{article}
\usepackage{geometry}
\usepackage[dvipdfmx]{graphicx}
\usepackage{amsfonts}
\usepackage{amsmath}
\usepackage{url}
\begin{document}
\title{Accelerated Multiple Precision Matrix Multiplication using Strassen's Algorithm and Winograd's Variant}
\author{Tomonori Kouya\thanks{Shizuoka Institute of Science and Technology}\\\url{http://na-inet.jp/}}
\date{2014-08-28}
%
%
\maketitle

\abstract{The Strassen algorithm and Winograd's variant accelerate matrix multiplication by using fewer arithmetic operations than standard matrix multiplication. Although many papers have been published to accelerate single- as well as double-precision matrix multiplication by using these algorithms, no research to date has been undertaken to accelerate multiple precision matrix multiplication. In this paper, we propose a multiple precision matrix multiplication program for matrices of any size and test its performance. We also reveal special properties of our program through its application to LU decomposition. }

%
\section{Introduction}

Current large-scale scientific computations use multiple precision (MP) floating-point arithmetic beyond the IEEE 754 single-precision (SP) and double-precision (DP) computation standard to obtain precise numerical solutions. Although MP arithmetic libraries, such as Multiple Precision Floating-Point Reliability (MPFR) and the GNU Multiple Precision Arithmetic Library (GMP), are software-based implementations, their MP numerical computations are typically much slower than hardware-based SP and DP computations. To prevent the consequent increase in computational cost, efficient MP numerical computation requires acceleration techniques, such as effective use of cache memory and algorithms to reduce the complexity of the computations.

Matrix multiplication is one of the most important parts of numerical computation. It is well known through research in DP matrix multiplication \cite{Huss-lederman96implementationof, strassen_cuda}, that its computational cost can be reduced by using Strassen's algorithm \cite{strassen_original} and Winograd's variant \cite{Coppersmith1990251}. By referring to past results, we can expect that MP matrix multiplication using these algorithms is more effective than in case of DP arithmetic. On the other hand, less precise numerical results may be obtained by applying Strassen's algorithm and its variant \cite{higham_accuracy}. 

In this paper, we propose the acceleration of MP matrix multiplication using Strassen's algorithm by comparing block matrix multiplication to increase the hit ratio of the cache memory in the CPU. We apply this accelerated MP matrix multiplication to LU decomposition, and examine both well-conditioned and ill-conditioned examples in order to study its numerical properties.

%
\section{Algorithms of Matrix Product}

We consider the real matrix multiplication $C := AB = [c_{ij}]$ $\in\mathbb{R}^{m\times n}$, where $A =[a_{ij}]$ $\in \mathbb{R}^{m\times l}$ and $B = [b_{ij}]$$\in \mathbb{R}^{l\times n}$ in this paper. We use the following algorithm to calculate $c_{ij}$:
\begin{equation}
	c_{ij} := \sum^l_{k = 1} a_{ik} b_{kj}. \label{eqn:matrix_mul_simple}
\end{equation}
Equation (\ref{eqn:matrix_mul_simple}) is called ``simple matrix multiplication" (``Simple," for short).

To increase the hit ratio of the cache memory in the processor, ``block matrix multiplication" (Block) with divided $A$ and $B$ are always used in well-tuned Basic Linear Algebra Subprogram (BLAS) libraries, such as the Automatically Tuned Linear Algebra Software (ATLAS) and the Intel Math Kernel. In this paper, we divide $A$ and $B$ into small $M L$ pieces of $A_{ik}$ and $L N$ pieces of $B_{kj}$, respectively. We can hence obtain blocked $C_{ij}$ by the following matrix multiplication:
\begin{equation}
	C_{ij} := \sum^L_{k = 1} A_{ik} B_{kj}. \label{eqn:matrix_mul_block}
\end{equation}

These simple and blocked matrix multiplication procedures have identical computational cost.

On the other hand, Strassen's algorithm to reduce the computational cost of matrix multiplication is recursive \cite{strassen_original}. For even-dimensional matrices $A$ and $B$ ($m$, $n$, and $l$ are even), we divide $A$ and $B$ as follows:

\begin{equation}
 A = \left[\begin{array}{cc}
	A_{11} & A_{12} \\
	A_{21} & A_{22}
\end{array}\right],\ B = \left[\begin{array}{cc}
	B_{11} & B_{12} \\
	B_{21} & B_{22}
\end{array}\right]. \label{eqn:ab_block}
\end{equation}

We calculate intermediate block matrices $P_i$ $(i = 1, 2, ..., 7)$ by using four divided $A_{ij}$s and $B_{ij}$s $(i, j = 1, 2)$ as follows:

\begin{eqnarray*}
	P_1 &:=& (A_{11} + A_{22}) (B_{11} + B_{22}) \\
	P_2 &:= &(A_{21} + A_{22}) B_{11} \\
	P_3 &:= &A_{11} (B_{12} - B_{22}) \\
	P_4 &:= &A_{22} (B_{21} - B_{11}) \\
	P_5 &:= &(A_{11} + A_{12}) B_{22} \\
	P_6 &:= &(A_{21} - A_{11}) (B_{11} + B_{12}) \\
	P_7 &:= &(A_{12} - A_{22}) (B_{21} + B_{22}).
\end{eqnarray*}

By using the values of $P_i$ above, we can calculate $C$ as blocked $C_{ij}$ $(i, j = 1, 2)$ as follows:

\begin{equation}
	C := \left[\begin{array}{cc}
		P_1 + P_4 - P_5 + P_7 & P_3 + P_5  \\
		P_2 + P_4 & P_1 + P_3 - P_2 + P_6
\end{array}\right] 
\end{equation}
By applying Strassen's algorithm to matrix multiplication, the number of real multiplication $\mbox{Mul}(m, l, n)$ and real addition-subtraction operations $\mbox{Addsub}(m, l, n)$ to calculate matrix $C$ using $A$ and $B$ is reduced as follows:
\begin{eqnarray*}
	\mbox{Mul}(m, l, n) &=& 7\mbox{Mul}(m/2, l/2, n/2)\\
	\mbox{Addsub}(m, l, n) &=& 5\mbox{Addsub}(m/2, l/2) \\
					     &+& 5\mbox{Addsub}(l/2, n/2) \\
					     &+& 8\mbox{Addsub}(m/2, n/2).
\end{eqnarray*}


Winograd proposed the self-titled ``Winograd's variant" (Winograd) algorithm that requires fewer matrix addition and subtraction operations than Strassen's algorithm\cite {Coppersmith1990251}. Winograd's variant is constructed with divided even-dimensional matrices in the same manner as in Strassen's algorithm (\ref{eqn:ab_block}), following which it computes matrix multiplication in the following three steps:

\begin{equation}
\begin{split}
	S_1 &:= A_{21} + A_{22},\ S_2 := S_1 - A_{11} \\
	S_3 &:= A_{11} - A_{21},\ S_4 := A_{12} - S_2 \\
	S_5 &:= B_{12} - B_{11},\ S_6 := B_{22} - S_5 \\
	S_7 &:= B_{22} - B_{12},\ S_8 := S_6 - B_{21}
\end{split} \label{eqn:winograd_mod1}
\end{equation}

\begin{equation}
\begin{split}
	M_1 &:= S_2 S_6,\ M_2 := A_{11} B_{11}, M_3 := A_{12} B_{21} \\
	M_4 &:= S_3 S_7,\ M_5 := S_1 S_5,\ M_6 := S_4 B_{22} \\
	M_7 &:= A_{22} S_8
\end{split} \label{eqn:winograd_mod2}
\end{equation}

\begin{equation}
	T_1 := M_1 + M_2,\ T_2 := T_1 + M_4 \label{eqn:winograd_mod3}
\end{equation}

Through (\ref{eqn:winograd_mod1})$\rightarrow$ (\ref{eqn:winograd_mod2}) $ \rightarrow$ (\ref{eqn:winograd_mod3}), we can obtain $C$ as follows:
\begin{equation}
	C := \left[\begin{array}{cc}
		M_2 + M_3 & T_1 + M_5 + M_6  \\
		T_2 - M_7 & T_2 + M_5
\end{array}\right] 
\end{equation}

Winograd's variant involves the following arithmetical operations:
\begin{eqnarray*}
	\mbox{Mul}(m, l, n) &=& 7\mbox{Mul}(m/2, l/2, n/2)  \\
	\mbox{Addsub}(m, l, n) &=& 4\mbox{Addsub}(m/2, l/2)  \\
	&+& 4\mbox{Addsub}(l/2, n/2) \\
	&+& 7\mbox{Addsub}(m/2, n/2).
\end{eqnarray*}

As we can observe, it can reduce a $\mbox{Addsub}(m/2, l/2)$, a $\mbox{Addsub}(l/2, n/2)$, and a $\mbox{Addsub}(m/2, n/2)$ operation.

In addition to Strassen's algorithm and Winograd's variant, we implement two matrix multiplication algorithms: a simple three-loop algorithm (\ref{eqn:matrix_mul_simple}) and a block algorithm (\ref{eqn:ab_block}). The four algorithms can obtain matrix products of any precision for matrices of any size. Strassen and Winograd recursively divided matrices $A$ and $B$ until the row and column dimensions were smaller than $n_{min}$ as the minimal dimension. In case of an odd number of row or column dimensions of $A$ or $B$, we fit them to become even by using a mixture of dynamic padding and peeling \cite{Huss-lederman96implementationof}.  

\tablename\ \ref{table:relative_complexity} shows the reduction rates of Strassen's algorithm and Winograd's variant in comparison with Simple and Block in case of $n_{min} = 32$. Both recursive algorithms can reduce multiplication operations by 45\% and addition-subtraction operations by 52\% in case of $m \times n = 2048 \times 2048$. There is a difference of a few percentage points in the efficiency of the addition and subtraction operations between the Strassen algorithm and Winograd's variant, and manifests itself as a more significant difference in computational time, as shown in the next section.

\begin{table}[htbp]
  \centering
  \caption{Relative complexity of Strassen's and Winograd's algorithms (vs. Simple and Block algorithms)}\label{table:relative_complexity}
    \begin{tabular}{r|cc|cc}
    \hline
       &  \multicolumn{2}{c|}{Strassen} & \multicolumn{2}{c}{Winograd} \\
    \hline
    $n_{min} = 32$&  Add \& Sub & Mul   & Add \& Sub & Mul   \\ \hline
      255 $\times$ 255 & 0.678 & 0.781 & 0.678 & 0.764 \\
      256 $\times$ 256 & 0.670 & 0.772 & 0.670 & 0.755 \\
      257 $\times$ 257 & 0.674 & 0.775 & 0.674 & 0.758 \\
      511 $\times$ 511 & 0.590 & 0.688 & 0.590 & 0.672 \\
      512 $\times$ 512 & 0.586 & 0.684 & 0.586 & 0.668 \\
      513 $\times$ 513 & 0.589 & 0.686 & 0.589 & 0.670 \\
     1023 $\times$ 1023 & 0.514 & 0.605 & 0.514 & 0.590 \\
     1024 $\times$ 1024 & 0.513 & 0.603 & 0.513 & 0.588 \\
     1025 $\times$ 1025 & 0.514 & 0.604 & 0.514 & 0.589 \\
     2047 $\times$ 2047 & 0.449 & 0.531 & 0.449 & 0.517 \\
     2048 $\times$ 2048 & 0.449 & 0.530 & 0.449 & 0.516 \\
     2049 $\times$ 2049 & 0.450 & 0.531 & 0.450 & 0.517 \\ \hline
    \end{tabular}%
  \label{tab:Relative Complexity of Simple, Block, Strassen's, and Winograd's alogorithms}%
\end{table}%

%
\section{Benchmark tests of square and rectangle matrix multiplications}

In this section, we use $a_{ij}$ and $b_{ij}$, the elements of $A$ and $B$, respectively, as follows:
\[ a_{ij} = \sqrt{5} (i + j - 1),\ b_{ij} = \sqrt{3} (n - i + 1).  \]
We then show the results of $C := AB$. Our numerical computational environment was as follows:
\begin{description}
\item[H/W] Intel Core i7 3850 (3.6 GHz), 64 GB RAM
\item[S/W] Scientific Linux 6.3 x86\_64, Intel C Compiler Ver. 13.0.1, BNCpack ver. 0.8, MPFR 3.1.2, GMP 5.1.3
\end{description}
All computations were serially executed without any parallelization. Since MPFR is a binary multiple precision floating-point library, we used a binary length of the mantissa as precision within a range of $128$ to $8192$ bits. 

We first discuss the results of square matrix multiplication ($m = n = l$). \tablename\ \ref{table:block} is obtained by using Simple and Block. Block($n_{min}$) represents the minimal dimension of the divided block matrices $A_{ik}$ and $B_{kj}$ as three values of $n_{min} = 16, 32, 64$. Block is the most effective algorithm in the case of 128-bit precision arithmetic, but we cannot recognize the difference between Simple and Block in case of 1024-bit precision. In case of 128 bits, the largest relative error in the elements of $C$ was $1.34 \times 10^{-37}$ and the smallest was $5.23\times 10^{-39}$. As a results, we obtained many times of smallest computational times in the case of $n_{min} = 32$.

\begin{table}[htbp]
  \centering
  \caption{Computation time: Block algorithm (128 bits)}\label{table:block}
    \begin{tabular}{r|rrrr}
    \hline
    \multicolumn{1}{c|}{$m\times n$} & Simple & Block(16) & Block(32) & Block(64) \\
    \hline
     255 $\times$ 255 & 1.06  & 1.20  & 1.22  & 1.24 \\
     256 $\times$ 256 & 1.25  & 1.22  & 1.22  & 1.25 \\
     257 $\times$ 257 & 1.04  & 1.25  & 1.28  & 1.37 \\
     511 $\times$ 511 & 9.60  & 9.71  & 9.61  & 10.02 \\
     512 $\times$ 512 & 10.83 & 9.70  & 9.68  & 9.96  \\
     513 $\times$ 513 & 10.02 & 9.89  & 9.97  & 10.44 \\
    1023 $\times$ 1023 & 107.78 & 77.63 & 77.80 & 79.36 \\
    1024 $\times$ 1024 & 213.09 & 77.77 & 77.72 & 79.51 \\
    1025 $\times$ 1025 & 94.62 & 78.92 & 78.41 & 81.48 \\
    2047 $\times$ 2047 & 756.81 & 627.75 & 619.21 & 648.31 \\
    2048 $\times$ 2048 & 1679.04 & 624.86 & 618.87 & 639.71 \\
    2049 $\times$ 2049 & 632.74 & 623.24 & 625.69 & 640.84 \\
    \hline
    \end{tabular}%
\end{table}%

We show \tablename\ \ref{table:strassen_winograd128bits} (128 bits precision) and \tablename\ \ref{table:strassen_winograd1024bits} (1024 bits) for the sake of comparison. These results of Strassen's algorithm and Winograd's variant are obtained with $n_{min} = 32$ due to these of Blocks.

\begin{table}[htbp]
  \centering
  \caption{Computation time: Strassen's and Winograd's algorithms (128 bits)}\label{table:strassen_winograd128bits}
    \begin{tabular}{r|rrr}
    \hline
    \multicolumn{1}{c|}{$m\times n$} & min(Simple, Block) & Strassen & Winograd \\
	\hline
     255 $\times$ 255 & 1.06  & 0.72  & 0.63 \\
     256 $\times$ 256 & 1.22  & 0.70  & 0.57 \\
     257 $\times$ 257 & 1.04  & 0.74  & 0.60 \\
     511 $\times$ 511 & 9.60  & 4.84  & 4.06 \\
     512 $\times$ 512 & 9.68  & 4.77  & 3.73 \\
     513 $\times$ 513 & 9.89  & 4.92  & 3.88 \\
    1023 $\times$ 1023 & 77.63 & 32.02 & 25.57 \\
    1024 $\times$ 1024 & 77.72 & 31.53 & 24.10 \\
    1025 $\times$ 1025 & 78.41 & 32.21 & 24.77 \\
    2047 $\times$ 2047 & 619.21 & 211.80 & 163.87 \\
    2048 $\times$ 2048 & 618.87 & 211.19 & 155.67 \\
    2049 $\times$ 2049 & 623.24 & 212.79 & 157.52 \\
    \hline
    \end{tabular}%
\end{table}%

\begin{table}[htbp]
  \centering
  \caption{Computation time: Strassen's and Winograd's algorithms (1024 bits)}\label{table:strassen_winograd1024bits}
    \begin{tabular}{r|rrr}
    \hline
    \multicolumn{1}{c|}{$n\times n$} & min(Simple, Block) & Strassen & Winograd \\
	\hline
     255 $\times$ 255 & 5.46  & 2.33  & 1.95 \\
     256 $\times$ 256 & 5.53  & 2.31  & 1.72 \\
     257 $\times$ 257 & 5.61  & 2.41  & 1.81 \\
     511 $\times$ 511 & 43.81 & 13.40 & 10.57 \\
     512 $\times$ 512 & 44.20 & 13.02 & 9.44 \\
     513 $\times$ 513 & 44.37 & 13.38 & 9.81 \\
    1023 $\times$ 1023 & 352.79 & 76.93 & 57.98 \\
    1024 $\times$ 1024 & 355.99 & 74.58 & 52.47 \\
    1025 $\times$ 1025 & 356.58 & 76.36 & 54.22 \\
    2047 $\times$ 2047 & 2820.16 & 454.02 & 329.41 \\
    2048 $\times$ 2048 & 2824.34 & 446.87 & 302.56 \\
    2049 $\times$ 2049 & 2829.95 & 456.08 & 307.05 \\
    \hline
    \end{tabular}%
\end{table}%

The maximum relative errors in $c_{ij}$s are as follows:
\begin{description}
	\item[128bits] Strassen: $3.20 \times 10^{-36}$, Winograd: $2.25 \times 10^{-35}$
	\item[1024bits] Strassen: $6.30 \times 10^{-306}$, Winograd: $3.92 \times 10^{-305}$
\end{description}
On occasion, the results of Winograd's variant are worse that those of Strassen's algorithm by one decimal digit.

We list the computation times of rectangle matrix multiplication for the four algorithms in \tablename\ \ref{table:rectangle}. All matrix multiplication operations obtain 1024-dimensional square matrices as their final result. In these cases, Block is faster than Simple beyond $l = 255$ or 511, and Winograd is always faster than Strassen within 37 seconds.

\begin{table}[htbp]
  \centering
  \caption{Computation time: Rectangle matrix multiplication ($n_{min} = 32$, Unit: seconds)}\label{table:rectangle}
\begin{tabular}{r|rrrr} \hline
  \multicolumn{5}{c}{128 bits computation}\\
 $m( = n), l $ & Simple & Block(32) & Strassen & Winograd \\ \hline
  1024,   63 & 5.09  & 5.86  & 5.31  & 4.29 \\
  1024,   64 & 5.18  & 5.91  & 5.03  & 3.99 \\
  1024,   65 & 5.26  & 6.30   & 5.24  & 4.21 \\
  1024,  127 & 10.47 & 11.8  & 9.43  & 6.45 \\
  1024,  128 & 10.54 & 11.86 & 8.71  & 5.42 \\
  1024,  129 & 10.63 & 12.24 & 8.90   & 5.65 \\
  1024,  255 & 52.51 & 23.69 & 14.67 & 9.39 \\
  1024,  256 & 52.27 & 23.77 & 13.13 & 7.51 \\
  1024,  257 & 52.56 & 24.07 & 13.39 & 7.70 \\
  1024,  511 & 110.75 & 47.42 & 26.40  & 16.55 \\
  1024,  512 & 106.19 & 47.44 & 21.85 & 11.83 \\
  1024,  513 & 110.82 & 47.96 & 22.08 & 12.09 \\ \hline
\multicolumn{5}{c}{ }  \\ \hline
\multicolumn{5}{c}{1024 bits computation}\\
    $m( = n), l$ &  Simple & Block(32) & Strassen & Winograd \\ \hline
    1024,  63 & 24.71 & 26.74 & 19.92 & 14.88 \\
    1024,  64 & 25.77 & 27.13 & 19.49 & 14.31 \\
    1024,  65 & 26.24 & 27.82 & 20.04 & 14.91 \\
    1024, 127 & 52.37 & 53.75 & 32.26 & 19.02 \\
    1024, 128 & 53.17 & 54.23 & 30.69 & 16.74 \\
    1024, 129 & 53.59 & 54.97 & 30.79 & 17.37 \\
    1024, 255 & 105.04 & 108.90 & 47.43 & 24.15 \\
    1024, 256 & 106.34 & 108.64 & 43.12 & 19.83 \\
    1024, 257 & 106.92 & 109.83 & 43.63 & 20.48 \\
    1024, 511 & 244.07 & 216.82 & 71.30 & 34.18 \\
    1024, 512 & 245.69 & 216.65 & 61.01 & 25.22 \\
    1024, 513 & 247.96 & 217.69 & 62.14 & 25.91 \\ \hline
    \end{tabular}%
\end{table}%

%
\section{Application to LU decomposition}

It is well known that matrix multiplication can be applied to LU decomposition \cite{golub4thed}. In this section, no LU decomposition involves any pivoting.

We consider the linear equation (\ref{eqn:linear_eq}) with $A \in\mathbb{R}^{n\times n}$, $\mathbf{b}\in\mathbb{R}^n$ 
\begin{equation}
	A\mathbf{x} = \mathbf{b}. \label{eqn:linear_eq}
\end{equation}
We use direct methods for the LU decomposition of the coefficient matrix by setting the block size to $K$. LU decomposition with matrix multiplication (the underlined part) is then as follows: 
\begin{enumerate}
	\item Divide $A$ into $A_{11}\in\mathbb{R}^{K\times K}$, $A_{12}\in\mathbb{R}^{K\times (n - K)}$, $A_{21}\in\mathbb{R}^{(n - K)\times K}$, and $A_{22}\in\mathbb{R}^{(n-K)\times (n-K)}$.
	\item Decompose $A_{11}$ into $L_{11} U_{11} (= A_{11})$, and then transform $A_{12}$ to $U_{12}$ and $A_{21}$ to $L_{21}$.
	\item $A^{(1)}_{22} := A_{22} - \underline{L_{21} U_{12}}$
\end{enumerate}
After substituting $A := A^{(1)}_{22}$, repeat the above algorithm until $n - K \geq 0$.

We employ a random matrix as an instance of a well-conditioned matrix and a Lotkin matrix as that of an ill-conditioned one. 
\begin{description}
	\item[Random Matrix] $a_{ij}$ is a random number in $[-1, 1]$.
	\item[Lotkin Matrix] $a_{ij} = \begin{cases}
		1 & (i = 1) \\
		1 / (i + j - 1) & (i \geq 2) \\
	\end{cases}$
\end{description}
The true solution is $\mathbf{x}$ $ = [0\ 1\ ...\ n-1]^T$, and we set $\mathbf{b} := A\mathbf{x}$. The condition numbers $\|A\|_1\|A^{-1}\|_1$ of the random matrix and the Lotkin matrix in $n = 1024$ are $4.4\times 10^{6}$ and $4.3\times 10^{1576}$, respectively.  For the Lotkin matrix, we must use more than 8192 bits (about 2466 decimal digits )  in $n = 1024$.

\begin{figure}[htb]
\begin{center}
\includegraphics[width=.8\textwidth]{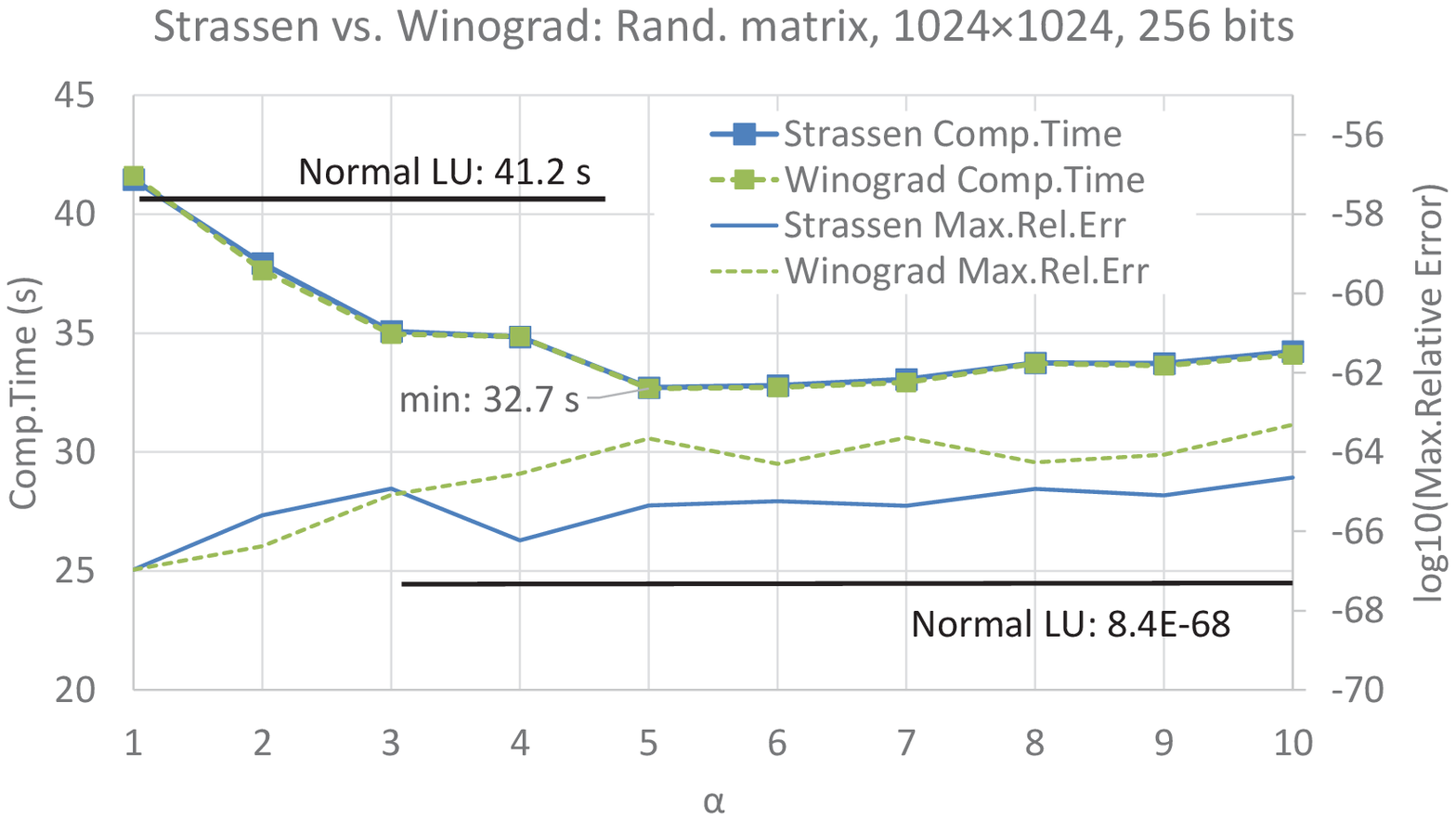}
\includegraphics[width=.8\textwidth]{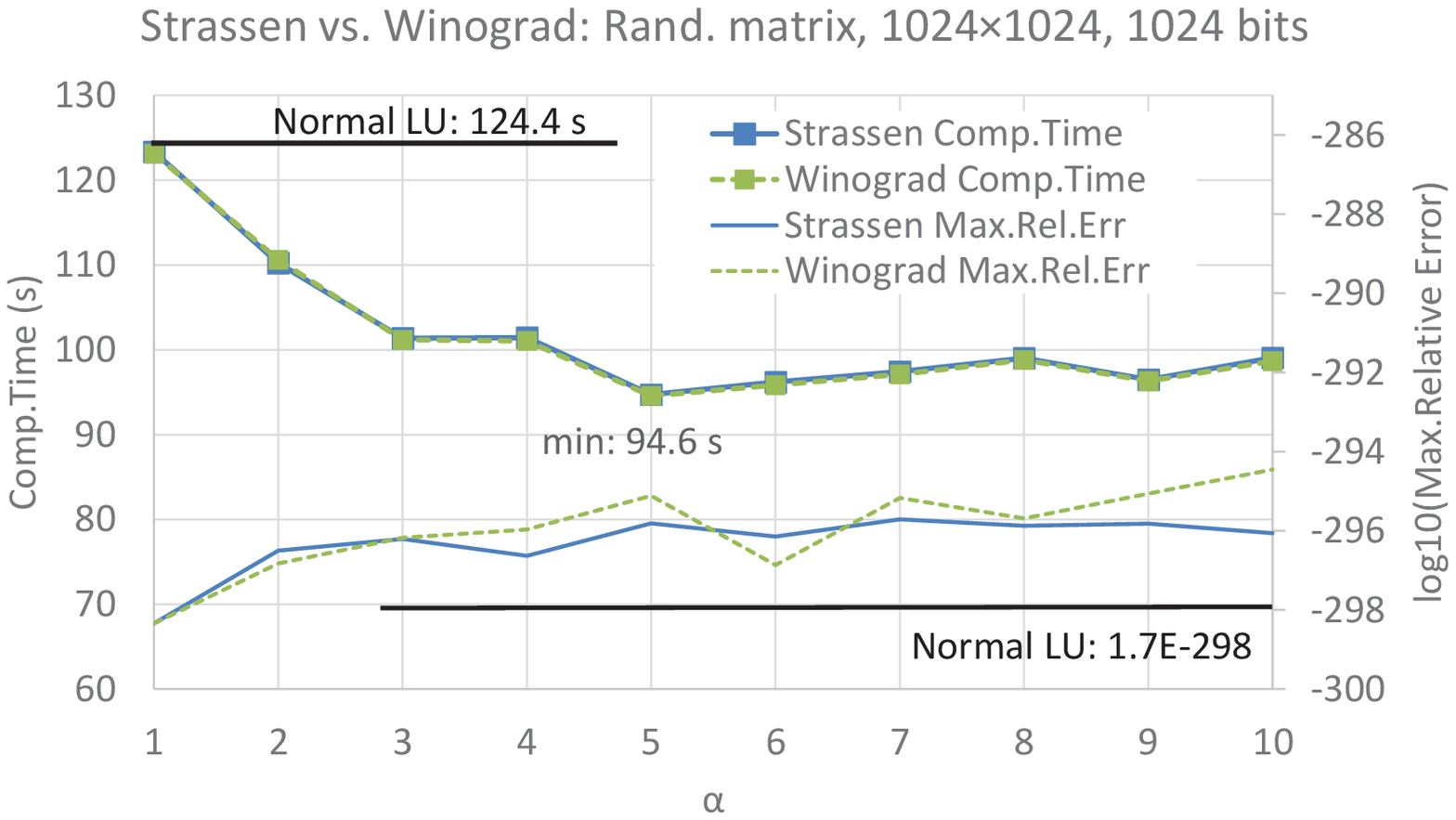}
\caption{Computation time and relative error of $1024 \times 1024$ random matrix (Upper: 256 bits, Lower: 1024 bits)}\label{fig:lu_strassen_random}
\end{center}
\end{figure}

\begin{figure}[htb]
\begin{center}
\includegraphics[width=.8\textwidth]{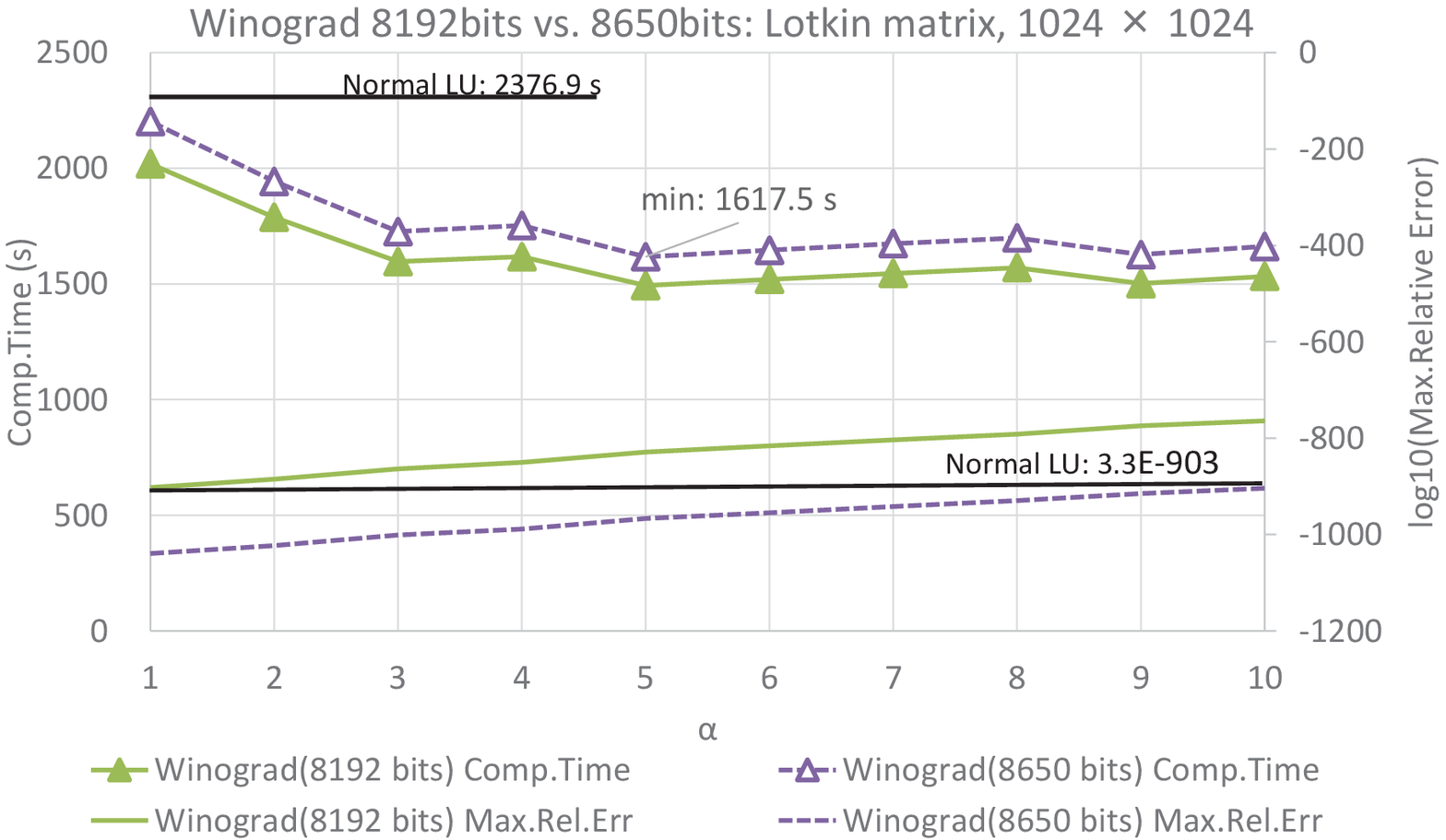}
\caption{Computation time and relative error of $1024 \times 1024$ Lotkin matrices (8192 and 8650 bits)}\label{fig:lu_strassen_lotkin}
\end{center}
\end{figure}

The size of the $K$s are set as $K = \alpha n_{min}$ ($\alpha = 1, 2, ..., 10$) and $n_{min} = 32$. Furthermore, we investigated the computation time (seconds) and the maximum relative error of the numerical solutions $\mathbf{x}$ at each $\alpha$. \figurename\ref{fig:lu_strassen_random} (random matrix) and \figurename\ref{fig:lu_strassen_lotkin} (Lotkin matrix) show the results. For comparison, the computation time and the maximum relative errors obtained using normal LU decomposition (column-wise LU) are shown in these figures.

We observe that we can reduce computation time by 21 to 26\% for a random matrix ($n = 1024$). For larger values of $\alpha$, the maximum relative errors grow from approximately two to four decimal digits. The computation times of Strassen's algorithm and Winograd's variant are within two seconds of each other.

We only show the results of using Winograd variant on the Lotkin matrix. In this case, the relative error increased 138 decimal digits ($n = 1024$). Thus, Winograd's variant operates in 8650 bits of computation in order to recover the increment of the relative error. Consequently we can reduce the computation time by 32 \%.

%
\section{Conclusion and future work}

We obtained the following results through our benchmark tests involving Simple, Block, Strassen's algorithm, and Winograd's variant. 
\begin{itemize}

\item The Block algorithm was more efficient than Simple algorithm when precision was relatively low, even in a multiple precision arithmetic environment.

\item Winograd's variant is always faster than Strassen's algorithm.

\item LU decomposition with Strassen's algorithm and Winograd's variant is faster than column-wise LU, but causes the loss of significant digits when the relevant coefficient matrix is ill-conditioned, such as Lotkin matrix.

\end{itemize}

In future research,  we will modify the block and recursive algorithms by using turning and parallelizing techniques.

%

\end{document}